\documentclass[12pt,a4paper]{amsart}
\usepackage{amsmath, amssymb, amsfonts,enumerate}
\usepackage[all]{xy}
\usepackage{amscd}

\newtheorem{theorem}{Theorem}[section]
\newtheorem{lemma}[theorem]{Lemma}
\newtheorem{corollary}[theorem]{Corollary}
\newtheorem{proposition}[theorem]{Proposition}

\theoremstyle{definition}

\newtheorem*{definition*}{Definition}

\theoremstyle{remark}

\newtheorem*{remarks}{Remarks}

\newtheorem*{examples}{Examples}

\numberwithin{equation}{section}

\newcommand {\K}{\mathbb{K}} 

\newcommand{\PP}{\mathcal{P}}
\newcommand{\CC}{\mathcal{C}}

\newcommand{\UU}{\mathcal{U}}
\newcommand{\NN}{\mathcal{N}}

\newcommand{\BB}{\mathcal{B}}

\DeclareMathOperator{\Fix}{Fix}
\DeclareMathOperator{\CA}{CA}
\DeclareMathOperator{\LCA}{LCA}
\DeclareMathOperator{\Ker}{Ker}

\begin{document}
\title{Expansive actions on uniform spaces and surjunctive maps}

\author{Tullio Ceccherini-Silberstein}
\address{Dipartimento di Ingegneria, Universit\`a del Sannio, C.so
Garibaldi 107, 82100 Benevento, Italy}
\email{tceccher@mat.uniroma1.it}
\author{Michel Coornaert}
\address{Institut de Recherche Math\'ematique Avanc\'ee,
UMR 7501,                                             Universit\'e  de Strasbourg et CNRS,
                                                 7 rue Ren\'e-Descartes,
                                               67000 Strasbourg, France}
\email{coornaert@math.unistra.fr}
\subjclass[2000]{54E15, 37B15, 68Q80}
\keywords{Expansive action, uniform space, surjunctive map, cellular automaton, shift action, subshift, group algebra, stable finiteness}
\date{January 11, 2011.}
\begin{abstract}
We present a uniform  version of a result of M. Gromov on the 
surjunctivity of maps commuting with expansive group actions
and discuss several applications. We prove in particular that for any
group $\Gamma$ and any field $\K$, the space of $\Gamma$-marked groups
$G$ such that the group algebra $\K[G]$ is stably finite is compact.
  \end{abstract}

\maketitle

\section{Introduction}

A map $f$ from a set $X$ into itself is said to be \emph{surjunctive} if it is surjective or not injective
 \cite{gromov}.
Thus, a non-surjunctive map is a map which is injective but not surjective.
With this terminology at hand, Dedekind's characterization of infinite sets may be rephrased by saying that a set $X$ is infinite if and only if it admits a non-surjunctive map $f \colon X \to X$.
 Similarly, elementary linear algebra tells us that a vector space is infinite-dimensional if and only if it admits a non-surjunctive endomorphism.
In fact, it turns out that the absence of non-surjunctive endomorphisms for a given mathematical object $X$ often reflects some ``finiteness" property of $X$.
\par
The word \emph{surjunctive}
was created by W. Gottschalk \cite{gottschalk} who introduced the notion of a surjunctive group. Let $G$ be a group. Given a set $A$, consider the set $A^G$ equipped with the \emph{prodiscrete} topology, that is, with the product topology obtained by taking the discrete topology on each factor $A$ of $A^G$. 
  There is a natural action of $G$  on $A^G$ 
  obtained by shifting coordinates via left multiplication in $G$ (see Section \ref{sec:ca}).
This action is called the $G$-\emph{shift} and the study of its dynamical properties is the central theme of symbolic dynamics.
The group $G$ is said to be \emph{surjunctive} if, for any finite set $A$, every $G$-equivariant continuous map $\tau \colon A^G \to A^G$ is surjunctive. 
Gottschalk asked whether every group $G$ is surjunctive. Although there is a large class of groups, namely sofic groups (see \cite{weiss}, \cite{elek}), which are known to be surjunctive, this question remains open.
Here, it is worth mentioning that the existence of a non-sofic group is also an open problem.  
\par
In algebraic geometry, there is a famous theorem of J. Ax \cite{ax} which says that every endomorphism of a complex algebraic variety is surjunctive.
Ax' surjunctivity theorem was extended by M. Gromov
 \cite{gromov} to certain classes of ``symbolic algebraic spaces"
 obtained by taking projective limits of algebraic varieties.
In order to establish these generalizations,
 Gromov developed various techniques for proving surjunctivity of maps 
$f \colon X \to X$ obtained as limits of surjunctive maps $f_i \colon X_i \to X_i$.
He also used one of these techniques to prove the closedness of the set of marked surjunctive groups (see \cite{gromov} and \cite{GG} for an alternative proof based on model theory).
The purpose of our paper is to give a detailed exposition of  this technique in a somewhat more general setting that the one considered in \cite{gromov}.
To be a little bit more precise, we extend results on the surjunctivity of maps commuting with expansive actions on metric spaces by replacing metric spaces by uniform spaces.
It turns out that this extension is quite natural and we think it may help
to a better understanding of the roles of the various mathematical concepts involved in Gromov's approach.
 \par
 Uniform spaces were introduced by A. Weil \cite{weil} as a generalization of metric spaces and topological groups. In a uniform space $X$, the closeness of a pair of points is not measured by a real number, like in a metric space, but by the fact that
 this pair of points belong or does not belong to certain subsets of the Cartesian square $X \times X$ which are called 
  the \emph{entourages} of the uniform structure 
 (see Section \ref{sec:uniform} for precise definitions).
 Every metric on a set $X$ defines a uniform structure on $X$ and every uniform structure defines a topology on the underlying set. In the case when the uniform structure comes from a metric, the associated
 topology coincides with the topology defined by the metric.
  
  Uniformly continuous maps between uniform spaces are defined as being the maps which pullback entourages to entourages. 
\par
There is a natural uniform structure on the set $\PP(X)$ of subsets of a uniform space $X$, which
is called the \emph{Hausdorff-Bourbaki} uniform structure. Its name comes from the fact that  it generalizes the Hausdorff distance on the set of bounded closed subsets of a metric space and was first mentioned in exercises in the treatise of N. Bourbaki \cite{bourbaki}.
\par
The classical definition of expansivity for groups actions on metric spaces may be extended to uniform spaces.
 More precisely, an action of a group $\Gamma$ on a uniform space $X$ is said to be \emph{expansive} if there is an entourage $V_0$ such that there is no pair of distinct points whose orbit under the diagonal action of $\Gamma$ on $X \times X$ remains inside $V_0$.
 \par 
 One of our main results is the following theorem whose proof will be given in Section \ref{sec:gil}:

\begin{theorem}
\label{t;A} 
Let $X$ be a uniform space equipped with a uniformly continuous and expansive action of a 
group $\Gamma$ and 
let $f \colon X \to X$ be a uniformly continuous and $\Gamma$-equivariant map.
Let $Y$ be a compact subset of $X$.
 Suppose that there exists a net $(Z_i)_{i \in I}$ of $\Gamma$-invariant subsets of $X$, which converges to $Y$ in the Hausdorff-Bourbaki topology such that $f(Z_i) \subset Z_i$ and the restriction maps $f\vert_{Z_i} \colon Z_i \to Z_i$ are surjunctive for all $i \in I$.   
Then $Y$ is $\Gamma$-invariant, $f(Y) \subset Y$, and the restriction map $f\vert_Y \colon Y \to Y$ is surjunctive.   
\end{theorem}

As an immediate consequence of the preceding theorem, we get

\begin{corollary}
\label{c:surj-subs-closed}
Let $X$ be a uniform space equipped with a uniformly continuous and expansive action of a 
group $\Gamma$ and 
let $f \colon X \to X$ be a uniformly continuous and $\Gamma$-equivariant map.
Let $\Sigma(f)$ denote the set of $\Gamma$-invariant compact subsets $Y \subset X$ such that $f(Y) \subset Y$
and the restriction map $f\vert_Y \colon Y \to Y$ is surjunctive.
Then $\Sigma(f)$ is closed in the set of all subsets of $X$ for the Hausdorff-Bourbaki topology.
\qed
\end{corollary}

Consider now a group $G$ and an arbitrary set $A$.
There is a natural uniform structure on $A^G$ which is called the \emph{prodiscrete} uniform structure (see Subsection \ref{ss;product}).
The associated topology is the prodiscrete topology on $A^G$ and the shift action is expansive and uniformly continuous. Note that if
$G$ is uncountable then the prodiscrete topology on $A^G$ is not metrizable as soon as $A$ has at least two elements.
The uniformly continuous $G$-equivariant maps $\tau \colon A^G \to A^G$ are precisely the \emph{cellular automata} over the group
$G$ and the alphabet $A$ (see Section \ref{sec:ca}). 

Let $\Gamma$ be a group. The set of quotients of $\Gamma$ may be identified with the set $\NN(\Gamma)$ of normal
subgroups of $\Gamma$. The set $\NN(\Gamma)$ is called the set of $\Gamma$-\emph{marked groups}. 
It is a closed (and hence compact) subset of the compact space $\PP(\Gamma) = \{0,1\}^\Gamma$ for the prodiscrete topology.

By applying
Theorem \ref{t;A} to the shift on $A^\Gamma$ with $A$ finite, we deduce in Section \ref{sec:mg} the above mentioned closedness theorem for
surjunctive groups:
\begin{corollary}[Gromov]
\label{c:surjunct-closed}
Let $\Gamma$ be a group. 
Then the set of normal subgroups $N \subset \Gamma$ such that the quotient group $\Gamma/N$ is surjunctive is closed (and hence compact) in $\NN(\Gamma)$.  
\end{corollary}

Let now $\K$ be a field.
A group $G$ is said to be L$_\K$-\emph{surjunctive} if, for any finite-dimensional vector space $V$ over $\K$, every linear cellular automaton
$\tau \colon V^G \to V^G$ is surjunctive.  It is known that every sofic group is L$_\K$-surjunctive (see \cite{CC-inj}).
In Section \ref{sec:lca}, we shall establish the closedness of the set of L$_\K$-surjunctive marked groups:
\begin{theorem}
\label{t:L-surjunct-closed}
Let $\Gamma$ be a group and let $\K$ be a field.
Then the set of normal subgroups $N \subset \Gamma$ such that the quotient group $\Gamma/N$ is L$_\K$-surjunctive is closed (and hence compact) in $\NN(\Gamma)$.  
\end{theorem}
The proof of the preceding theorem uses a stronger version of Theorem \ref{t;A} 
(see Corollary \ref{c;gil}), where the compactness hypothesis for $Y$ is replaced by the weaker hypotheses that both $Y$ and $f(Y)$ are closed in $X$ together with the fact that the restriction of $f$ to $Y$ is a uniform embedding.
\par
One of the famous conjectures of I. Kaplansky \cite{kaplansky} about the algebraic properties of group algebras is that $\K[G]$ is stably finite
for any group $G$ and any field $\K$ (we recall that a ring $R$ is said to be \emph{stably finite} if every one-sided invertible square matrix over $R$ is also two-sided invertible). This conjecture remains open up to now although it was settled by Kaplansky himself for all fields of characteristic $0$ and later by Elek and Szab\'o \cite{ES} (see also \cite{CC-inj}) for sofic groups in any characteristic.
In \cite{CC-inj}, it is shown that a group $G$ is L$_\K$-surjunctive if and only if the group algebra $\K[G]$ is stably finite. We thus have:
\begin{corollary}
Let $\Gamma$ be a group and let $\K$ be a field.
Then the set of normal subgroups $N \subset \Gamma$ such that the group algebra $\K[\Gamma/N]$ is stably finite is closed (and hence compact) in $\NN(\Gamma)$.  
\end{corollary}

The paper is organized as follows.
 Section \ref{sec:uniform} contains background material on uniform spaces.
The definition of the Hausdorff-Bourbaki uniform structure on the set of subsets of a uniform space is recalled in Section \ref{sec:hausdorff-bourbaki}.
For the convenience of the reader, proofs of some well known facts about the Hausdorff-Bourbaki uniform structure that are needed in the sequel are also included in this section. In Section \ref{sec:gil} we prove a uniform version of Lemma 4.H" in \cite{gromov} which is used to prove 
Theorem \ref{t;A}. Basic properties of cellular automata are presented in Section \ref{sec:ca}. 
Section \ref{sec:mg} is devoted to marked groups. We show that if $G = \Gamma/N$ is a $\Gamma$-marked group, then the set $A^G$ is canonically
isomorphic, as a uniform $G$-space, to the subset $\Fix(N) \subset A^\Gamma$ consisting of the points of $A^\Gamma$ which are fixed by the
normal subgroup $N \subset \Gamma$. We also prove that the map $N \mapsto \Fix(N)$ yields a uniform embedding of $\NN(\Gamma)$, the space of $\Gamma$-marked groups, into $\PP(A^\Gamma)$. This gives a symbolic representation of the space of $\Gamma$-marked groups which is used, in combination
with Theorem \ref{t;A} (resp. Corollary \ref{c:surj-subs-closed}), in the
proof of Corollary \ref{c:surjunct-closed} (resp. Theorem \ref{t:L-surjunct-closed}) which is completed in Section \ref{sec:mg} 
(resp. Section \ref{sec:lca}). 

\section{Uniform structures}
\label{sec:uniform}

In this section we collect some basic facts about uniform spaces, uniformly continuous maps, and uniformly continuous group actions. 
More details and proofs may be found for example in \cite{bourbaki}, \cite{james}, or \cite{kelley}. 

 Throughout the paper, we shall use the following general notation.
 Let $X$ be a set. 
 We denote by $\PP(X)$ the set of all subsets of $X$ 
and by $\Delta_X$ the diagonal in $X \times X$, that is,
$\Delta_X = \{ (x,x) : x \in X \}$.
Let $R$ be a binary relation on $X$, that is,
 a subset $R \subset X \times X$. Given an element $y \in X$, we define the subset $R[y] \subset X$ by
$R[y] = \{x \in X : (x,y) \in R \}$. For $Y \subset X$, we define
the set $R[Y] \subset X$ by $R[Y] = \bigcup_{y \in Y} R[y]$, so that we have in particular $R[y] = R[\{y\}]$ for all $y \in X$.
\par
The \emph{inverse} $\overset{-1}{R}$ of $R$ is the binary relation
$\overset{-1}{R} = \{ (x,y) : (y,x) \in R \}$.
One says that $R$ is \emph{symmetric} if $\overset{-1}{R} = R$.
\par  
 The \emph{composite} of two binary relations $R$ and $S$ is the binary relation
  $$
R \circ S = \{ (x,y): \text{ there exists  } z \in X \text{  such that  } (x,z) \in R \text{  and  } (z,y) \in S \}  \subset X \times X.
$$

\subsection{Uniform spaces}
\label{ss:us}

\begin{definition*}
Let $X$ be a set. A \emph{uniform structure}\index{uniform ! --- structure} on $X$ is a non--empty set
$\UU$ of subsets
of $X \times X$ satisfying the following conditions:
\begin{enumerate}[(UN-1)]
\item if $V \in \UU$, then $\Delta_X \subset V$;
\item if $V \in \UU$ and $V \subset V' \subset X \times X$, then $V' \in \UU$;
\item if $V \in \UU$ and $W \in \UU$, then $V \cap W \in \UU$;
\item if $V \in \UU$, then $\overset{-1}{V} \in \UU$;
\item if $V \in \UU$, then there exists $W \in \UU$ such that $W
\circ W \subset V$.
\end{enumerate}
\end{definition*}

 A set $X$ equipped with a uniform structure $\UU$ is called a {\it 
uniform space} and
the elements of $\UU$ are called the \emph{entourages}  of $X$.
\par
Let $(X,\UU)$ be a uniform space.
If $Y$ is a subset of $X$, then $\UU_Y = \{V \cap (Y \times Y) : V \in \UU \}$ is a uniform structure on $Y$, which is said to be 
\emph{induced} by $\UU$.
\par
A subset
$\BB \subset \UU$ is called a \emph{base} of $\UU$ if
for each $V \in \UU$ there exists
$B \in \BB$ such that $B \subset V$.
\par
It is easy to see that a base $\BB$ satisfies the 
following properties:
\begin{enumerate}[\rm (B-1)]
\item 
if $V \in \BB$, then $\Delta_X \subset V$;
\item 
if $V \in \BB$ and $W \in \BB$, then there exists $U \in \BB$ 
such that $U \subset V \cap W$;
\item 
if $V \in \BB$, then there exists $W \in \BB$ such that $W 
\subset \overset{-1}{V}$;
\item 
if $V \in \BB$, then there exists $W \in \BB$ such that $W 
\circ W \subset V$.
\end{enumerate}
Conversely, given a set $X$ and a subset $\BB \subset \PP(X \times X)$ satisfying conditions (B-1) -- (B-4), there exists a unique uniform structure $\UU$ on
$X$ admitting $\BB$ as a base.

There is a topology on $X$ which is associated with the uniform structure. 
This is the topology for which the neighborhoods of an arbitrary point $x \in X$ consist of the sets $V[x]$, where $V$ runs over all entourages of $X$.
This topology is Hausdorff if and only if the intersection of the entourages of $X$ coincides with the diagonal $\Delta_X \subset X \times X$.

\begin{examples}
\label{examples:discrete}
1) The \emph{discrete uniform structure}
 on a set $X$ is the uniform 
structure whose entourages
consist of all subsets of $X \times X$ containing $\Delta_X$.
The topology associated with this uniform structure is the discrete topology on $X$.
\par
2) If $(X,d)$ is a metric space, then there is a uniform structure on $X$ associated with the metric. A base for this uniform structure is given by the sets
$V_\varepsilon = \{(x,y) \in X \times X : d(x,y) < \varepsilon\}$, $\varepsilon > 0$. The topology associated with this uniform structure coincides with the topology associated with the metric.
\end{examples}

\subsection{Uniformly continuous maps}
\label{ss:uniformly-continuous}
Let $X$ and $Y$ be uniform spaces.
\par 
A map $f \colon X \to Y$ is said to be
\emph{uniformly continuous}  if it satisfies the following condition: 
for each entourage $W$
of $Y$, there exists an entourage $V$ of $X$ such that $(f \times f)(V) \subset W$. Here $f \times f$ denotes the map from $X \times X$ into $Y \times Y$ defined by
$(f \times f)(x_1,x_2) = (f(x_1),f(x_2))$ for all $(x_1,x_2) \in X \times X$.
Every uniformly continuous map $f \colon X \to Y$ is continuous but the converse fails to hold in general.
\par
A \emph{uniform isomorphism} between uniform spaces $X$ and $Y$ is a bijective map $f \colon X \to Y$ such that both $f$ and $f^{-1}$ are uniformly continuous.
One says that a map $f \colon X \to Y$ is a \emph{uniform embedding} if $f$ is injective and induces a uniform isomorphism between $X$ and $f(X)$.

\subsection{Product of uniform spaces}
\label{ss;product}
Suppose that $(X_\lambda)_{\lambda \in \Lambda}$ is a family of uniform spaces.
Then the smallest uniform structure on the Cartesian product $X = \prod_{\lambda \in \Lambda} X_\lambda$ for which all projection maps $\pi_\lambda \colon X \to X_\lambda$, $\lambda \in \Lambda$, are uniformly continuous
is called the \emph{product uniform structure} on $X$.
A base of entourages for the product uniform structure on $X$ is obtained by taking all subsets of 
$X \times X$ which are of the form
\begin{align*}
\prod_{\lambda \in \Lambda} V_\lambda 
&\subset \prod_{\lambda \in \Lambda} X_\lambda \times X_\lambda \\
&= \left( \prod_{\lambda \in \Lambda} X_\lambda \right) \times \left( \prod_{\lambda \in \Lambda} X_\lambda \right) \\
&= X \times X,
\end{align*}
where $V_\lambda \subset X_\lambda \times X_\lambda$ is an entourage of $X_\lambda$ and $V_\lambda = X_\lambda \times X_\lambda$ for all but finitely many $\lambda \in \Lambda$.
\par
If $Y$ is a uniform space, then a map $f \colon Y \to X$ is uniformly continuous if and only if the maps $\pi_\lambda \circ f \colon Y \to X_\lambda$ are uniformly continuous for all $\lambda \in \Lambda$.
\par
In the particular case when each $X_\lambda$ is endowed with the discrete uniform structure, the product uniform structure on $X = \prod_{\lambda \in \Lambda} X_\lambda$ is called the
\emph{prodiscrete uniform structure} on $X$.

\subsection{Uniformly continuous and expansive actions}
Suppose that $X$ is a uniform space equipped with an action of a group $\Gamma$.
\par
 One says that the action of $\Gamma$ on $X$ is
\emph{uniformly continuous} if the map $f_\gamma \colon X \to X$ defined by $f_\gamma(x) = \gamma x$ is uniformly continuous for each $\gamma \in \Gamma$.
\par
Consider the diagonal action of $\Gamma$ on $X \times X$ defined by $\gamma (x,y) = (\gamma x, \gamma y)$ for all $x,y \in X$ and $\gamma \in \Gamma$.
One says that the action of $\Gamma$ on $X$ is \emph{expansive} if there exists an entourage $W_0$ of $X$ satisfying the following property:
for any two distinct points $x,y \in X$, there exists an element $\gamma \in \Gamma$ such that $\gamma(x, y) \notin W_0$. In other words, the action of $\Gamma$ on $X$ is expansive if and only if there exists an entourage $W_0$ such that
$$
 \bigcap_{\gamma \in \Gamma} \gamma^{-1} W_0 = \Delta_X.
$$
Note that if $X$ admits a uniformly continuous and expansive action then the topology on $X$ is necessarily Hausdorff.

\section{The Hausdorff-Bourbaki uniform structure}
\label{sec:hausdorff-bourbaki}

In this section we briefly review the definition and basic properties of the Hausdorff-Bourbaki uniform structure on the set of subsets of a uniform space 
(see \cite[ch. II exerc. 5 p. 34 and exerc. 6 p. 36]{bourbaki}).
Some proofs of well known facts that will be used in the next section are included for the convenience of the reader.
 \par
Let $X$ be a uniform space with uniform structure $\UU$.
 For $V \in \UU$, we define the subset $\widehat{V} \subset \PP(X) \times \PP(X)$ by
\begin{equation}
\label{e:V-hat}
\widehat{V} = \left\{(Y,Z) \in \PP(X) \times \PP(X) : Z \subset V[Y] \text{  and  } Y \subset V[Z] \right\}.
\end{equation}
One easily checks that the set $\BB = \{\widehat{V} : V \in \UU \} \subset \PP(\PP(X) \times \PP(X))$ satisfies conditions (B-1) -- (B-4). Therefore, there exists a unique uniform 
structure on $\PP(X)$ admitting $\BB$ as a base.
This uniform structure is called the \emph{Hausdorff-Bourbaki uniform structure} associated with $\UU$.
The topology associated with the Hausdorff-Bourbaki uniform structure is called the \emph{Hausdorff-Bourbaki topology} on $\PP(X)$.
If $(Y_i)_{i \in I}$ is a net in $\PP(X)$ (that is, a family of subsets of $X$ indexed by some directed set $I$) and $Z \in \PP(X)$, we shall write
$Y_i \overset{H-B}{\to} Z$ to mean that the net $(Y_i)_{i \in I}$ converges to $Z$ in the Hausdorff-Bourbaki topology. 

\begin{proposition}
\label{p:grom-haus-sep-closed}
Let $X$ be a uniform space and let $Y$ and $Z$ be closed subsets of $X$.
Suppose that there is a net $(T_i)_{i \in I}$ of subsets of $X$ such that
$T_i \overset{H-B}{\to} Y$ and $T_i \overset{H-B}{\to} Z$.
 Then one has $Y = Z$.
\end{proposition}

\begin{proof}
Let $y \in Y$ and let $\Omega$ be a neighborhood of $y$ in $X$.
Then there is a symmetric entourage $V$ of $X$ such that $V[y] \subset \Omega$.
Choose an entourage $W$ of $X$ such that $W \circ W \subset V$.
Since the net $(T_i)_{i \in I}$ converges to both $Y$ and $Z$, we can find an element $i_0 \in I$ such that $Y \subset W[T_{i_0}]$ and $T_{i_0} \subset W[Z]$.
Thus, there exist $t \in T_{i_0}$ and $z \in Z$ such that
$(y,t) \in W$ and $(t,z) \in W$. This implies
$(y,z) \in W \circ W \subset V$ and hence $(z,y) \in V$ since $V$ is symmetric. 
  It follows that $z \in V[y] \subset \Omega$. 
This shows that $y$ is in the closure of $Z$.
Since $Z$ is closed in $X$, we deduce that $Y \subset Z$. By symmetry, we also have $Z \subset Y$. Consequently, $Y = Z$. 
\end{proof}

\begin{remarks}
1) An immediate consequence of Proposition \ref{p:grom-haus-sep-closed} is that, if $X$ is a uniform space, then the topology induced by the Hausdorff-Bourbaki topology on the set of closed subsets of $X$ is Hausdorff.
\par
2) Suppose that $(X,d)$ is a metric space and let $\CC_b(X)$ denote the subset of $\PP(X)$ consisting of all closed bounded subsets of $X$. For $x \in X$ and $r > 0$, denote by $B(x,r)$ the open ball of radius $r$ centered at $x$.
The \emph{Hausdorff metric} is the metric $\delta$ on $\CC_b(X)$ defined by
 $$
\delta(Y,Z) = \inf\{ r > 0 : Z \subset \bigcup_{y \in Y} B(y,r) \text{  and  } Y \subset \bigcup_{z \in Z} B(z,r) \}.
$$
 The uniform structure on $\CC_b(X)$ associated with the metric $\delta$ is the uniform structure induced by the Hausdorff-Bourbaki structure on $\PP(X)$. 
 \end{remarks}

\begin{proposition}
Let $X$ be a uniform space and let us equip $\PP(X)$ with the Hausdorff-Bourbaki uniform structure and $\PP(X) \times \PP(X)$ with the corresponding product uniform structure.
Then the map $\eta \colon \PP(X) \times \PP(X) \to \PP(X)$ defined by $\eta(Y,Z) = Y \cup Z$ is uniformly continuous.
\end{proposition}

\begin{proof}
Let $V$ be an entourage of $X$ and consider the entourage $\widehat{V}$ of $\PP(X)$ given by \eqref{e:V-hat}.
 Suppose that $Y_1, Y_2,Z_1, Z_2$ are subsets of $X$ such that $(Y_1,Y_2) \in \widehat{V}$ and $(Z_1,Z_2) \in \widehat{V}$. Then we have $Y_1 \subset V[Y_2] \subset V[Y_2 \cup Z_2]$ and $Z_1 \subset V[Z_2] \subset V[Y_2 \cup Z_2]$, and
therefore $Y_1 \cup Z_1 \subset V[Y_2 \cup Z_2]$. Similarly, we get $Y_2 \cup Z_2 \subset V[Y_1 \cup Z_1]$. We deduce that $(Y_1 \cup Z_1,Y_2 \cup Z_2) \in \widehat{V}$.
Consequently, we have $(\eta \times \eta)(W) \subset \widehat{V}$,
where $W$ is the entourage of $\PP(X) \times \PP(X)$ defined by
$$
W = \{ ((Y_1,Z_1),(Y_2,Z_2)) \in (\PP(X) \times \PP(X)) \times (\PP(X) \times \PP(X)) : (Y_i,Z_i) \in \widehat{V} (i = 1,2) \}. 
$$
This shows that $\eta$ is uniformly continuous.
\end {proof}

\begin{corollary}
\label{c:limit-inclusion-gh}
Let $X$ be a uniform space. Let $(Y_i)_{i \in I}$ and $(Z_i)_{i \in I}$ be nets in $\PP(X)$ such that
$Y_i \subset Z_i$ for all $i \in I$. Suppose that there exist closed subsets $A$ and $B$ of $X$ such that $Y_i \overset{H-B}{\to} A$ and $Z_i \overset{H-B}{\to} B$. Then one has $A \subset B$.
\end{corollary}

\begin{proof}
We have $Y_i \cup Z_i \overset{H-B}{\to} A \cup B$ by continuity of the union map $\eta$.
Since $Y_i \cup Z_i = Z_i$ for all $i \in I$, Proposition \ref{p:grom-haus-sep-closed}
implies that $A \cup B = B$, 
that is, $A \subset B$.
\end{proof}

 \begin{proposition}
\label{p:uc-implies-gh-uc}
Let $X$ and $Y$ be uniform spaces and let $f \colon X \to Y$ be a uniformly continuous map.
Then the map
$f_* \colon \PP(X) \to \PP(Y)$ which sends each subset $A \subset X$ to its image $f(A) \subset Y$ is uniformly continuous with respect to the Hausdorff-Bourbaki uniform structures on $\PP(X)$ 
and $\PP(Y)$.
\end{proposition}

\begin{proof}
Let $W$ be an entourage of $Y$ and let
$$
\widehat{W} = \{ (B_1,B_2) \in \PP(Y) \times \PP(Y) : B_2 \subset W[B_1] \text{  and  } B_1 \subset W[B_2] \}
$$
be the associated entourage of $\PP(Y)$. Since $f$ is uniformly continuous, there is an entourage $V$ of $X$ such that
$$
(f \times f)(V) \subset W.
$$
Suppose that $(A_1,A_2) \in \widehat{V}$, that is, $A_2 \subset V[A_1]$ and $A_1 \subset V[A_2]$.
If $a_1 \in A_1$, then there exists $a_2 \in A_2$ such that $(a_1,a_2) \in V$ and hence $(f(a_1),f(a_2)) \in W$. Therefore, we have $f_*(A_1) \subset W[f_*(A_2)]$. Similarly, we get $f_*(A_2) \subset W[f_*(A_1)]$. 
It follows that $(f_*(A_1),f_*(A_2)) \in \widehat{W}$.
This shows that
$(f_* \times f_*)(\widehat{V}) \subset \widehat{W}$. Consequently, $f_*$ is uniformly continuous.   
\end{proof}

 \begin{corollary}
\label{c:limit-f-inv}
Let $X$ be a uniform space and let $f \colon X \to X$ be a uniformly continuous map.
Let $Y$ be a subset of $X$ such that the sets $Y$ and $f(Y)$ are both closed in $X$.
Suppose that there is a net $(Z_i)_{i \in I}$ of subsets of $X$ such that $f(Z_i) \subset Z_i$ for all $i \in I$  and $Z_i \overset{H-B}{\to} Y$.
Then one has $f(Y) \subset Y$.  
\end{corollary}

\begin{proof}
We have $f(Z_i) \subset Z_i$ for all $i \in I$. 
Since $f(Z_i) \overset{H-B}{\to} f(Y)$ by Proposition \ref{p:uc-implies-gh-uc}, we deduce that $f(Y) \subset Y$ by applying Corollary \ref{c:limit-inclusion-gh}.
\end{proof}

In particular, we have:

\begin{corollary}
\label{c:limit-Gamma-inv}
Let $X$ be a uniform space equipped with a uniformly continuous action of a group $\Gamma$.
Let $Y$ be a closed subset of $X$.
Suppose that there is a net $(Z_i)_{i \in I}$ of $\Gamma$-invariant subsets of $X$ such that $Z_i \overset{H-B}{\to} Y$.
Then $Y$ is $\Gamma$-invariant.
\qed
 \end{corollary}

\section{Gromov's injectivity lemma}
\label{sec:gil}

The following result is a uniform version of Lemma 4.H'' in \cite{gromov}. 

\begin{theorem}[Gromov's injectivity lemma]
\label{t;gil} 
Let $X$ be a uniform space equipped with a uniformly continuous and expansive action of a 
group $\Gamma$ and 
let $f \colon X \to X$ be a uniformly continuous and $\Gamma$-equivariant map.
Suppose that $Y$ is a subset of $X$ such that
the restriction of $f$ to $Y$ is a uniform embedding.
Then there exists an entourage $V$ of $X$ satisfying the following 
property: if $Z$ is a
$\Gamma$-invariant subset of $X$ such that $Z \subset V[Y]$, then
the restriction of $f$ to $Z$ is injective.
\end{theorem}

 \begin{proof}
By expansivity of the action of $\Gamma$ on $X$, there exists an entourage $W_0$ 
of $X$ such that
\begin{equation}  
\label{e:expansive2}
\bigcap_{\gamma \in \Gamma} \gamma^{-1} W_0 = \Delta_X.
\end{equation}
It follows from the axioms of a uniform structure that we can find a
symmetric entourage $S$ of $X$ such that
\begin{equation} 
\label{e:S-W}
S \circ S \circ S \subset W_0.
\end{equation}
Since the restriction of $f$ to $Y$ is a uniform
embedding, there exists an entourage $T$ of $X$ such that 
\begin{equation} 
\label{e:T-S}
(f(y_1),f(y_2)) \in T \Rightarrow (y_1,y_2) \in S
\end{equation}
for all $y_1,y_2 \in Y$.
Let $U$ be a symmetric entourage of $X$ such that
\begin{equation} \label{e:U-T}
U \circ U \subset T.
\end{equation}
As $f$ is uniformly continuous, we can find an entourage $E$ of
$X$ such that  
\begin{equation} 
\label{e:E-U}
(x_1,x_2) \in E \Rightarrow (f(x_1),f(x_2)) \in U
\end{equation}
for all $x_1,x_2 \in X$.
\par
Let us show that  the entourage $V = S \cap E$ has the required property.
Suppose that $Z$ is a 
$\Gamma$-invariant subset of $X$ such that $Z \subset V[Y]$.  Let us show that the restriction of $f$ to 
$Z$ is injective.
\par
Let $z'$ and $z''$ be points in $Z$ such that $f(z') = f(z'')$.
Since $f$ is $\Gamma$-equivariant, we have
\begin{equation}
\label{e;z',zz''}
  f(\gamma z') = f(\gamma z'')
\end{equation} 
for all $\gamma \in \Gamma$. 
As $\gamma z'$ and $\gamma z''$ stay in $Z$ and $Z \subset V[Y]$, we can find points
$y_\gamma'$ and $y_\gamma''$ in $Y$ such that $(\gamma z',y_\gamma') \in V$
and $(\gamma z'',y_\gamma'') \in V$.
Since $V \subset E$, it follows from \eqref{e:E-U} that
$(f(\gamma z'),f(y_\gamma'))$ and $(f(\gamma z''),f(y_\gamma''))$ are
in $U$. As $U$ is symmetric we also have $(f(y'_\gamma), f(\gamma z')) \in U$ and therefore
$(f(y_\gamma'),f(y_\gamma'')) \in U \circ U \subset T$ by
\eqref{e;z',zz''} and \eqref{e:U-T}. By applying \eqref{e:T-S}, we deduce that
$(y_\gamma',y_\gamma'') \in S$. On the other hand, we also have 
$(\gamma z', y'_\gamma) \in S$
and $(y_\gamma'', \gamma z'') \in S$ since $V \subset S$ and $S$ is symmetric. It follows that
$$(\gamma z',\gamma z'')  \in S \circ S \circ S \subset W_0
$$
by \eqref{e:S-W}.
This gives us
$$
(z',z'') \in \bigcap_{\gamma \in \Gamma} \gamma^{-1} W_0,
$$
and hence $z' = z''$ by \eqref{e:expansive2}. This shows that the
restriction of $f$ to $Z$ is injective.
\end{proof}

Let $f \colon X \to X$ be a map from a set $X$ into itself. 
If $Y$ is a subset of $X$ such that $f(Y) \subset Y$,  we denote by $f\vert_Y \colon Y \to Y$ the restriction map given by $f\vert_Y(y) = f(y)$ 
for all $y \in Y$.

\begin{corollary}
\label{c;gil} 
Let $X$ be a uniform space equipped with a uniformly continuous and expansive action of a 
group $\Gamma$ and 
let $f \colon X \to X$ be a uniformly continuous and $\Gamma$-equivariant map.
Let $Y$ be a subset of $X$ such that
 $Y$ and $f(Y)$ are both closed in $X$ and the restriction of $f$ to $Y$ is a uniform embedding.
Suppose that there exists a net $(Z_i)_{i \in I}$ of $\Gamma$-invariant subsets of $X$ such that $f(Z_i) \subset Z_i$  and the restriction maps $f\vert_{Z_i} \colon Z_i \to Z_i$ are surjunctive for all $i \in I$ and
$Z_i \overset{H-B}{\to} Y$. 
Then $Y$ is $\Gamma$-invariant and one has $f(Y) = Y$.  
 \end{corollary}

\begin{proof}
The $\Gamma$-invariance of $Y$ directly follows from Corollary \ref{c:limit-Gamma-inv}.
Let $V$ be an entourage of $X$ as in Theorem \ref{t;gil}.
As $Z_i \overset{H-B}{\to} Y$, there is an element $i_0 \in I$ such that $Z_i \subset V[Y]$ for all $i \geq i_0$. Using the fact that $f\vert_{Z_i}$ is surjunctive,
we deduce that $f(Z_i) = Z_i$ for all $i \geq i_0$. 
We have $f(Z_i) \overset{H-B}{\to} f(Y)$ by Proposition \ref{p:uc-implies-gh-uc}.
As $Y$ and $f(Y)$ are closed in $X$,
we conclude that $f(Y) = Y$ by applying Proposition \ref{p:grom-haus-sep-closed}. 
\end{proof}

\begin{proof}[Proof of Theorem \ref{t;A}]
First observe that the space $X$ is Hausdorff since it admits an expansive uniformly continuous action. It follows that both
$Y$ and $f(Y)$ are closed in $X$ by compactness of $Y$.

 The fact that $Y$ is $\Gamma$-invariant and satisfies $f(Y) \subset Y$ is a consequence of
 Corollary \ref{c:limit-Gamma-inv} and Corollary \ref{c:limit-f-inv}.
 \par
 Suppose now that $f\vert_Y$ is injective.
As $Y$ is compact, we deduce that the restriction of $f$ is a uniform embedding. 
We conclude that $f\vert_Y$ is surjective by applying Corollary \ref{c;gil}.
\end{proof}

\section{Cellular automata}
\label{sec:ca}

Let $G$ be a group and let $A$ be an arbitrary set that we shall call the \emph{alphabet}. We consider the set $A^G$ consisting of all maps $x \colon G \to A$.
The elements of $A^G$ are called \emph{configurations} over the group $G$ and the alphabet $A$.
\par
Given a subset $S \subset G$, 
 we denote by $\pi_S \colon A^G \to A^S$  the projection map. 
 \par
We equip the set $A^G = \prod_{g \in G} A$ with its prodiscrete uniform structure (cf. Subsection \ref{ss;product}) and the left action 
of $G$
 defined by
$gx(h) = x(g^{-1}h)$ for all $g,h \in G$ and $x \in A^G$.   This action is called the $G$-\emph{shift} on $A^G$.  
\begin{proposition} 
\label{p;UCE}
 The $G$-shift on $A^G$ is
uniformly continuous and expansive.
 \end{proposition}

\begin{proof}
Uniform continuity follows from the fact that $G$ acts on $A^G$ by 
permuting coordinates. Expansiveness is due to the fact that the 
action of $G$ on itself by left multiplication
is transitive. Indeed, consider the entourage $W_0$ of $A^G$ defined by
$$W_0 = \{(x,y) \in A^G \times A^G: x(1_G) = y(1_G)\}.$$
Given $g \in G$, we have $(x,y) \in g^{-1}W_0$ if and only if $x(g^{-1}) = 
y(g^{-1})$. Thus $\bigcap_{g \in G} g^{-1}W_0$
is equal to the diagonal in $A^G \times A^G$.
\end{proof}

The space $A^G$ with the $G$-shift action is called the \emph{full shift} on $G$ and $A$.  A $G$-invariant subset $X \subset A^G$ is called a \emph{subshift} (note that here we do not require $X$ to be closed in $A^G$).

\begin{definition*}
Let $X, Y \subset A^G$ be two subshifts.
A \emph{cellular automaton}  from $X$ to $Y$ is a map $\tau \colon X \to Y$ 
satisfying the following property:
there exist a finite subset $S \subset G$ 
and a map $\mu \colon \pi_S(X) \to A$ such that 
\begin{equation} 
\label{def:automate}
 \tau(x)(g) = \mu(\pi_S(g^{-1}x)) \quad \text{for all } x \in X \text{ and } g \in G.
\end{equation}
Such a set $S$ is then called a \emph{memory set}
and $\mu$ is called a \emph{local defining map}
for $\tau$.
\end{definition*}

The following statement is proved in \cite{cs-c} (see also \cite{CC-book}) in the particular case when $X = Y = A^G$:

\begin{proposition} 
\label{p;CA}
Let $G$ be a group and let $A$ be a set. 
Let $X, Y \subset A^G$ be two subshifts and let $\tau \colon X \to Y$ be a map. 
Then the following conditions are equivalent:
\begin{enumerate}[\rm (a)]
\item 
$\tau$ is a cellular automaton;
\item 
$\tau$ is uniformly continuous and $G$-equivariant;
\item
There exists a cellular automaton $\widetilde{\tau} \colon A^G \to A^G$ such that $\widetilde{\tau}(X) \subset Y$ and $\tau(x) = \widetilde{\tau}(x)$
for all $x \in X$.
\end{enumerate}
\end{proposition}

\begin{proof}
Suppose that $\tau$ is a cellular automaton with memory set $S$ and local defining map $\mu \colon \pi_S(X) \to A$. Let $\widetilde{\mu} \colon A^S \to A$ be a map extending $\mu$.
Then the map $\widetilde{\tau} \colon A^G \to A^G$ defined by $\widetilde{\tau}(x)(g) = \widetilde{\mu}(\pi_S(g^{-1}x))$ for all $x \in A^G$ and $g \in G$ is a cellular automaton over $A^G$ such that
$\widetilde{\tau}(X) \subset Y$ and $\widetilde{\tau}(x) = \tau(x)$ for all $x \in X$. This shows that (a) implies (c).
\par
It is shown in \cite[Th. 1.1]{cs-c} that every cellular automaton $\widetilde{\tau} \colon A^G \to A^G$ is $G$-equivariant and uniformly continuous. Therefore (c) implies (b).
\par  
Suppose that $\tau \colon X \to Y$ is $G$-equivariant and uniformly continuous.
Since the map $f \colon X \to A$ defined by $f(x) = \tau(x)(1_G)$ is uniformly continuous, this implies that there is a finite subset $S \subset G$ such that if 
two elements $x,y \in X$ coincide on $S$ then $f(x) = f(y)$.
In other words, there is a map $\mu \colon \pi_S(X) \to A$ such that
$\tau(x)(1_G) = \mu(\pi_S(x))$ for all $x \in X$. Now the $G$-equivariance of $\tau$ gives us
$$
\tau(x)(g) = g^{-1}\tau(x)(1_G) = \tau(g^{-1}x)(1_G) = \mu(\pi_S(g^{-1}x))
$$ 
for all $x \in X$ and $g \in G$. This shows that (b) implies (a).
\end{proof}

If $X,Y \subset A^G$ are two subshifts, we shall denote by $\CA(X,Y;G,A)$ the set consisting of all cellular automata $\tau \colon X \to Y$.
If $Y = X$ we shall simply write $\CA(X;G,A)$ instead of $\CA(X,X;G,A)$.
Note that it immediately follows from the equivalence of conditions (a) and (b) in Proposition \ref{p;CA} that $\CA(X;G,A)$ is a monoid for the composition of maps.

\section{Marked groups}
\label{sec:mg}

Let $\Gamma$ be a group. 
A $\Gamma$-\emph{quotient} is a pair $(G,\rho)$,
 where
$G$ is a group and $\rho \colon \Gamma \to G$ is a group 
epimorphism. 
Two $\Gamma$-quotients   $(G_1,\rho_1)$ and $(G_2,\rho_2)$ 
are said to be \emph{equivalent} if there exists a group isomorphism $\phi \colon G_2 \to G_1$
such that the following diagram is commutative:
$$
 \xymatrix{ 
  & G_2 \ar[dd]^{\phi \cong} \\
 \Gamma \ar@{->>}[ru]^{\rho_2} \ar@{->>}[rd]_{\rho_1} & \\
 & G_1
 } 
 $$
 that is, such that $\rho_1 = \phi \circ \rho_2$.
An equivalence class of $\Gamma$-quotients is called a $\Gamma$-\emph{marked group}.
 Observe that two $\Gamma$-quotients 
 $(G_1,\rho_1)$ and $(G_2,\rho_2)$ are equivalent if and only if $\Ker(\rho_1) = \Ker(\rho_2)$.
Thus the set of $\Gamma$-marked groups may be identified with the set $\NN(\Gamma)$ consisting of all normal subgroups  of $\Gamma$.

Let us equip the set $\PP(\Gamma) = \{0,1\}^{\Gamma}$ with its prodiscrete uniform structure and $\NN(\Gamma)$ with the induced uniform structure.
A base of entourages of $\NN(\Gamma)$ is provided by the sets
$$
V_F = \{(N_1,N_2) \in \NN(\Gamma) \times \NN(\Gamma) : N_1 \cap F = N_2 \cap F \},
$$
where $F$ runs over all finite subsets of $\Gamma$.
Intuitively, two normal subgroups of $\Gamma$ are ``close" in $\NN(\Gamma)$ when their intersection with a large finite subset of $\Gamma$ coincide.

The space $\PP(\Gamma)$ is Hausdorff and totally disconnected since it is a product of discrete spaces.
Moreover, $\PP(\Gamma)$ is compact by the Tychonoff product theorem.
  The set $\NN(\Gamma)$  is closed in $\PP(\Gamma)$ and therefore $\NN(\Gamma)$ is a totally disconnected compact Hausdorff space (see \cite{champetier}, \cite{cornulier} and the references therein).   

Let $A$ be a set. Consider the set $A^\Gamma$ equipped with its prodiscrete uniform structure and the $\Gamma$-shift action.
For each $N \in \NN(\Gamma)$, let 
$$
\Fix(N) = \{ x \in A^\Gamma : \gamma x = x \text{ for all } \gamma \in N \} \subset A^\Gamma
$$
denote the set of points in $A^\Gamma$ which are fixed by $N$. 
Observe that $\Fix(N)$ is a closed subshift of $A^\Gamma$. Since $N$ acts trivially on $\Fix(N)$, the $\Gamma$-shift on $A^\Gamma$ induces an action of the quotient group
$G = \Gamma/N$ on $\Fix(N)$. 

\begin{proposition}
\label{p:describes-fix-N}
Let $N$ be a normal subgroup of a group $\Gamma$ and let $\rho \colon \Gamma \to G = \Gamma/N$ denote the canonical epimorphism.
Then $y \circ \rho \in \Fix(N)$ for every $y \in A^G$.
Moreover, the map $\rho^* \colon A^G \to \Fix(N)$ defined by $\rho^*(y) = y \circ \rho$ is 
a $G$-equivariant uniform isomorphism.
 \end{proposition}

 \begin{proof}
An element $x \in A^\Gamma$ is in $\Fix(N)$ if and only if $x(\nu^{-1}\gamma) = x(\gamma)$ for all $\gamma \in \Gamma$ and $\nu \in N$, that is, if and only if the configuration $x$ is constant on each coset of $\Gamma$ modulo $N$.
This proves the first assertion and the fact that the map $\rho^*$ is surjective. 
The injectivity of $\rho^*$ follows from the surjectivity of $\rho$.
Let $g \in G$ and choose $\gamma \in \Gamma$ such that $\rho(\gamma) = g$.
For all $y \in A^G$ and $\alpha \in \Gamma$, we have
\begin{align*}
g(y \circ \rho)(\alpha) &= \gamma (y \circ \rho)(\alpha) = y \circ \rho (\gamma^{-1}\alpha) = y(\rho(\gamma^{-1}\alpha)) = y(g^{-1}\rho(\alpha))\\
&= ((gy) \circ \rho) (\alpha).
\end{align*}
This shows that $g(y \circ \rho) = (gy) \circ \rho$, that is, $g\rho^*(y) = \rho^*(gy)$. 
Thus $\rho^*$ is $G$-equivariant.
\par
For each $\gamma \in \Gamma$, let $\pi_\gamma \colon A^\Gamma \to A$ and $\pi_\gamma' \colon A^G \to A$ denote the projection maps given by
$x \mapsto x(\gamma)$ and $y \mapsto y(\rho(\gamma))$ respectively.
The fact that $\pi_\gamma \circ \rho^* = \pi_\gamma'$ is uniformly continuous for all $\gamma \in \Gamma$ implies that the map $\rho^*$ is uniformly continuous.
Similarly, the uniform continuity of $(\rho^*)^{-1}$ 
follows from the fact that
$\pi_\gamma' \circ (\rho^*)^{-1} = \pi_\gamma\vert_{\Fix(N)}$ is uniformly continuous for each $\gamma \in \Gamma$. Consequently,
$\rho^*$ is a uniform isomorphism.  
  \end{proof}

By using the characterization of cellular automata provided by the equivalence of conditions (a) and (b) in Proposition \ref{p;CA}, we immediately deduce the following: 

\begin{corollary}
\label{c:rho-conj-ca}
If $\tau \colon A^G \to A^G$ is a cellular automaton over $A^G$, then the map $\tau^* \colon \Fix(N) \to \Fix(N)$ given by $\tau^* = \rho^* \circ \tau \circ (\rho^*)^{-1}$ is a cellular automaton over the subshift
$\Fix(N) \subset A^\Gamma$.
Moreover, the map $\Theta \colon \CA(A^G,G,A) \to \CA(\Fix(N),\Gamma,A)$ defined 
by $\Theta(\tau) = \tau^*$ is a monoid isomorphism.
\qed   
\end{corollary}

$$
\begin{CD}
A^G  @>\rho^*>> \Fix(N) \subset A^\Gamma \\
@V\tau VV  @VV{\tau^*}V \\
A^G @>>\rho^*> \Fix(N)
\end{CD}
$$

Let us equip $\PP(A^\Gamma)$ with the Hausdorff-Bourbaki uniform structure associated with the prodiscrete uniform structure on $A^\Gamma$.

\begin{theorem}
\label{t:psi-uc-n-p}
Let $\Gamma$ be a group and let $A$ be a set.
Then the map $\Psi \colon \NN(\Gamma) \to \PP(A^\Gamma)$ defined by $\Psi(N) = \Fix(N)$ is uniformly continuous.
Moreover, if $A$ contains at least two elements then $\Psi$ is a uniform embedding.
\end{theorem}

\begin{proof}
Let $N_0 \in \NN(\Gamma)$ and let $W$ be an entourage of $\PP(A^\Gamma)$.
Let us show that there exists an entourage $V$ of $\NN(\Gamma)$ such that
\begin{equation}
\label{e:unif-cont-psi-mg}
\Psi(V[N_0]) \subset W[\Psi(N_0)].
\end{equation}
This will prove that $\Psi$ is continuous.
By definition of the Hausdorff-Bourbaki uniform structure on $\PP(A^\Gamma)$, there is an entourage $T$ of $A^\Gamma$ such that
\begin{equation}
\label{e:ent-in-gr-haus}
\widehat{T} = \{ (X,Y) \in \PP(A^\Gamma) \times \PP(A^\Gamma) : Y \subset T[X] \text{ and } X \subset T[Y] \} \subset W.
\end{equation}
 Since $A^\Gamma$ is endowed with its prodiscrete uniform structure,
there is a finite subset $F \subset \Gamma$ such that
\begin{equation}
\label{e:u-f-t-entour}
U = \{ (x,y) \in A^\Gamma \times A^\Gamma : \pi_F(x) = \pi_F(y) \} \subset T,
\end{equation}
where $\pi_F \colon A^\Gamma \to A^F$ is the projection map.
Consider now the finite subset $E \subset \Gamma$ defined by
$$
E = F \cdot F^{-1} = \{ \gamma \eta^{-1} : \gamma,\eta \in F \},
$$
and the entourage $V$ of $\NN(\Gamma)$ given by
$$
V = \{ (N_1,N_2) \in \NN(\Gamma) \times \NN(\Gamma) : N_1 \cap E = N_2 \cap E \}.
$$
We claim that $V$ satisfies \eqref{e:unif-cont-psi-mg}. 
 To prove our claim, suppose that $N \in V[N_0]$. Let $x \in \Fix(N)$. The fact that $N \cap E = N_0 \cap E$ implies that if $\gamma$ and $\eta$ are elements of $F$ with $\gamma = \nu \eta$ for some $\nu \in N$, then $\nu \in N_0$ and therefore $x(\gamma) = x(\eta)$. Denoting by 
$\rho_0 \colon \Gamma \to \Gamma/N_0$ the canonical epimorphism, we deduce that we may find an element $x_0 \in A^{\Gamma/N_0}$ such that
$x(\gamma) = x_0 \circ \rho_0(\gamma)$ for all $\gamma \in F$.
We have $x_0 \circ \rho_0 \in \Fix(N_0)$ and $(x,x_0 \circ \rho_0) \in U$. Since $U \subset T$ by \eqref{e:u-f-t-entour}, this shows that $\Fix(N) \subset T[\Fix(N_0)]$.
Therefore $\Psi$ is continuous and hence uniformly continuous by compactness of $\NN(\Gamma)$. 
 \par
Suppose now that $A$ has at least two elements.
Let us show that $\Psi$ is injective.
  Let $N_1, N_2 \in \NN(\Gamma)$.
  Fix two elements $a,b \in A$ with $a \not= b$ and
consider the map $x \colon \Gamma \to A$ defined by $x(\gamma) = a$ if $\gamma \in N_1$ 
and $x(\gamma) = b$ otherwise. We clearly have $x \in \Fix(N_1)$.
Suppose that $\Psi(N_1) = \Psi(N_2)$, that is, $\Fix(N_1) = \Fix(N_2)$. Then for all $\nu \in N_2$, we have $\nu^{-1}x = x$ since $x \in \Fix(N_1) = \Fix(N_2)$, and hence $x(\nu) = \nu^{-1}x(1_\Gamma) = x(1_\Gamma) = a$. This implies $N_2 \subset N_1$. By symmetry, we also have $N_1 \subset N_2$.
Therefore $N_1 = N_2$. This shows that $\Psi$ is injective.
As $\NN(\Gamma)$ is compact   and $\PP(A^\Gamma)$ is Hausdorff, we conclude that $\Psi$ is a uniform embedding.     
 \end{proof}

\begin{proof}[Proof of Corollary \ref{c:surjunct-closed}]
Let $N \in \NN(\Gamma)$ and let $(N_i)_{i \in I}$ be a net in $\NN(\Gamma)$ converging to $N$.
Suppose that the groups $G_i = \Gamma/N_i$ are surjunctive for all $i \in I$ and let us show that the group $G = \Gamma/N$ is surjunctive as well.
\par
Let $A$ be a finite set and let $\tau \colon A^G \to A^ G$ be a cellular automaton over $A^G$.
Let $\rho \colon \Gamma \to G$ denote the canonical epimorphism and let 
    $\rho^* \colon A^G \to \Fix(N)$ be the map defined by
$\rho^* = y \circ \rho$ for all $y \in A^G$ (cf Proposition \ref{p:describes-fix-N}). 
By Corollary \ref{c:rho-conj-ca}, the map $\tau^* \colon \Fix(N) \to \Fix(N)$ given by $\tau^* = \rho^* \circ \tau \circ (\rho^*)^{-1}$ is a cellular automaton over the subshift $\Fix(N) \subset A^\Gamma$. 
By applying Proposition \ref{p;CA}, 
we deduce the existence of a cellular automaton $\widetilde{\tau} \colon A^\Gamma \to A^\Gamma$ such that $\widetilde{\tau}(\Fix(N)) \subset \Fix(N)$ and $\widetilde{\tau}\vert_{\Fix(N)} = \tau^*$.
\par
We claim that the hypotheses of Theorem \ref{t;A} are satisfied by taking $X = A^\Gamma$, $f = \widetilde{\tau}$, $Y = \Fix(N)$, and $Z_i = \Fix(N_i)$. First observe that the space $X$ is compact by the Tychonoff product theorem since $A$ is finite. As $Y$ is closed in $X$, it follows that $Y$ is compact. Moreover, the action of $\Gamma$ on $X$ is uniformly continuous and expansive by Proposition \ref{p;UCE}.
On the other hand, the cellular automaton  $f \colon X \to X$ is uniformly continuous and 
$\Gamma$-equivariant by Proposition \ref{p;CA}. 
We have $f(Z_i) \subset Z_i$ since the cellular automaton $f$ is $\Gamma$-equivariant. Moreover, $Z_i = \Fix(N_i)$ is $\Gamma$-invariant by normality of the subgroup $N_i \subset \Gamma$. Also, by Theorem \ref{t:psi-uc-n-p}, $Z_i \overset{H-B}{\to} Y$ since $N_i$ converges to $N$. 
For $i \in I$, let $\rho_i \colon \Gamma \to G_i = \Gamma/N_i$ denote the canonical epimorphism and let $\rho_i^* \colon A^{G_i} \to \Fix(N_i)$ be the
uniform isomorphism given by Proposition \ref{p:describes-fix-N}. 
It follows from Corollary \ref{c:rho-conj-ca} that the map $\tau_i = (\rho_i^*)^{-1} \circ \widetilde{\tau} \vert_{Z_i} \circ \rho_i^* \colon A^{G_i} \to A^{G_i}$ is a cellular automaton and hence is surjunctive since $G_i$ is a surjunctive group. 
This implies that $f\vert_{Z_i} = \widetilde{\tau}\vert_{\Fix(N_i)}$ is surjunctive as well.  
By applying Theorem \ref{t;A}, we deduce that $f\vert_Y = \tau^*$ and hence $\tau = (\rho^*)^{-1} \circ \tau^* \circ \rho^*$ is surjunctive. 
Consequently, the group $G$ is surjunctive.
    \end{proof}

 \section{Linear cellular automata}
\label{sec:lca}

Let $G$ be a group and  let $V$ be a vector space over a field $\K$.
The set $V^G$ has a natural vector space structure and the $G$-shift action on $V^G$ is linear (here and in the sequel ``linear" stands for ``$\K$-linear").
\par
A \emph{linear subshift} of $V^G$ is a $G$-invariant vector subspace of $V^G$.
\par
Let $X, Y \subset V^G$ be two linear subshifts. A \emph{linear cellular automaton} from $X$ to $Y$ is a cellular automaton $\tau \colon X \to Y$ which is linear.
Note that if $\tau \colon X \to Y$ is a cellular automaton with memory set $S \subset G$ and local defining map $\mu \colon \pi_S(X) \to V$, then
$\tau$ is linear if and only if $\mu$ is linear.

\begin{proposition} 
\label{p;LCA}
Let $G$ be a group and let $V$ be a vector space over a field $\K$. 
Let $X,Y \subset V^G$ be two linear subshifts and let $\tau \colon X \to Y$ be a map. 
Then the following conditions are equivalent:
\begin{enumerate}[\rm (a)]
\item 
$\tau$ is a linear cellular automaton;
\item 
$\tau$ is linear, uniformly continuous (with respect to the
uniform structures on $X$ and $Y$ induced by the prodiscrete uniform structure on $V^G$), 
and $G$-equivariant;
\item
$\tau$ is linear, continuous (with respect to the 
topologies on $X$ and $Y$ induced by the prodiscrete topology on $V^G$), and $G$-equivariant;
\item
$\tau$ is linear, continuous (with respect to the 
topologies on $X$ and $Y$ induced by the prodiscrete topology on $V^G$) at the constant configuration $x = 0$, and $G$-equivariant;
\item
there exists a linear cellular automaton $\widetilde{\tau} \colon V^G \to V^G$ such that $\widetilde{\tau}(X) \subset Y$ and $\tau(x) = \widetilde{\tau}(x)$ for all $x\in X$.
\end{enumerate}
\end{proposition}

\begin{proof}
The equivalences (a) $\Leftrightarrow$ (b) $\Leftrightarrow$ (e) immediately follow from Proposition \ref{p;CA} after observing that if $S \subset G$ then any linear map
$\mu \colon \pi_S(X) \to V$ may be extended to a linear map $\widetilde{\mu} \colon V^S \to V$. Since the topology associated with the prodiscrete uniform structure is the
prodiscrete topology (cf. Example (1) in Section \ref{ss:us})
and every uniformly continuous map is continuous
(cf. Section \ref{ss:uniformly-continuous}), we also have (b) $\Rightarrow$ (c).
The implication  (c) $\Rightarrow$ (d) is trivial. 
Let us show that (d) $\Rightarrow$ (a). Suppose that $\tau$ is linear, continuous at $0$, and
$G$-equivariant.
Then, the map $X \to V$ defined by $x \mapsto \tau(x)(1_G)$ is continuous at $0$ since the projection maps
$X \to V$ defined by $x \mapsto x(g)$ are continuous for all $g \in G$ and the composition of continuous maps is continuous.
We deduce that there exists a finite subset $S \subset G$  such that if $x \in X$  satisfies
$x(s) = 0$ for all $s \in  S$, then $\tau(x)(1_G) = 0$.
By linearity, we have that if two configurations $x$ and $y$ coincide on $S$ then $\tau(x)(1_G) = \tau(y)(1_G)$. Thus there exists a linear map $\mu \colon \pi_S(X) \to V$ such that $\tau(x)(1_G) = \mu(\pi_S(x))$ for all $x \in X$. As $\tau$ is $G$-equivariant, we deduce that $\tau(x)(g) = \tau(g^{-1}x)(1_G) = \mu(\pi_S(g^{-1}x))$ for all $x \in X$ and $g \in G$.
This shows that $\tau$ is the (linear) cellular automaton with memory set $S$ and local defining map $\mu$. Thus (d) implies (a). 
\end{proof}

We shall denote by $\LCA(X,Y;G,V)$ the set consisting 
of all linear cellular automata $\tau \colon X \to Y$. If $X = Y$ we shall simply write $\LCA(X;G,V)$
instead of $\LCA(X,X;G,V)$.
Note that $\LCA(X;G,V)$ has a natural structure of $\K$-algebra, where the multiplicative law is given by composition.
\par
If $N$ is a normal subgroup of a group $\Gamma$ and $V$ is a vector space over a filed $\K$, then the set $\Fix(N)$ of
points of $V^\Gamma$ which are fixed by $N$ is a closed linear subshift of $V^\Gamma$. From Proposition \ref{p:describes-fix-N}, we immediately deduce the following:

\begin{proposition}
\label{p:describes-fix-N-lin}
Let $N$ be a normal subgroup of a group $\Gamma$ and let $\rho \colon \Gamma \to G = \Gamma/N$ denote the canonical epimorphism.
Let $V$ be a vector space over a field $\K$.
Then $y \circ \rho \in \Fix(N)$ for every $y \in V^G$.
Moreover, the map $\rho^* \colon V^G \to \Fix(N)$ defined by $\rho^*(y) = y \circ \rho$ is 
a linear $G$-equivariant uniform isomorphism.
\qed
 \end{proposition}

\begin{corollary}
\label{c:rho-conj-lca}
If $\tau \colon V^G \to V^G$ is a linear cellular automaton over $V^G$, then the map $\tau^* \colon \Fix(N) \to \Fix(N)$ given by $\tau^* = \rho^* \circ \tau \circ (\rho^*)^{-1}$ is a linear cellular automaton over the linear subshift
$\Fix(N) \subset V^\Gamma$.
Moreover, the map $\Theta \colon \LCA(V^G,G,V) \to \LCA(\Fix(N),\Gamma,V)$ defined 
by $\Theta(\tau) = \tau^*$ 
is an isomorphism of $\K$-algebras.
\qed   
\end{corollary}

In the proof of Theorem \ref{t:L-surjunct-closed} we shall use the following results.
\begin{lemma}[Closure lemma] 
\label{l;closure}
Let $G$ be a group and let $V$ be a finite dimensional vector space over a field $\K$. Let
$\tau \colon V^G \to V^G$ be a linear cellular automaton. Then $\tau(V^G)$ is closed in $V^G$ for the prodiscrete topology.
\end{lemma}
\begin{proof} When $G$ is countable this is shown in \cite[Lemma 3.1]{CC-eden}. The general case is treated in 
\cite[Corollary 1.6]{CC-ind}.
\end{proof}

\begin{lemma}
\label{l;inverse}
Let $G$ be a group and let $V$ be a finite dimensional vector space over a field $\K$. Let
$\tau \colon V^G \to V^G$ be an injective linear cellular automaton and let $Y = \tau(V^G)$ denote the image of $\tau$. 
Then the inverse map $\tau^{-1}
\colon Y \to V^G$ is a linear cellular automaton and $\tau$ is
a uniform embedding.
\end{lemma}
\begin{proof} When $G$ is countable this is shown in \cite[Theorem 3.1]{CC-inj}. 
\par
We now treat the general case. Let $S \subset G$ be a memory set for $\tau$ and denote by $\mu \colon V^S \to V$ the corresponding
local defining map. Let $H \subset G$ be the subgroup generated by $S$. Note that $H$ is countable. Let $G/H$ denote the set of
all left cosets of $H$ in $G$. For $c \in G/H$ denote by 
$$
\pi_c \colon V^G \to V^c = \prod_{g \in c} V
$$ 
the projection map.
We have $V^G = \prod_{c \in G/H} V^c$ and, for every $x \in V^G$ we write $x = (x_c)_{c \in G/H}$, where $x_c = \pi_c(x) \in V^c$.

For $c \in G/H$, $g \in c$ and $z \in V^c$ we denote by $g^{-1}z \in V^H$ the configuration defined by
$$
(g^{-1}z)(h) = z(gh)
$$
for all $h \in H$. Then, if $x = (x_c)_{c \in G/H} \in V^G$ and $c \in G/H$, we have
$$
(g^{-1}x_c)(h) = x_c(gh) = x(gh) = (g^{-1}x)(h) = (g^{-1}x)_H(h)
$$
for all $h \in H$, that is,
\begin{equation}
\label{e;g-1-x-c}
g^{-1}x_c = (g^{-1}x)_H.
\end{equation}

For $c \in G/H$, define the map $\tau_c \colon V^c \to V^c$ by setting
$$
\tau_c(z)(g) = \mu(\pi_S^H(g^{-1}z))
$$
for all $z \in V^c$ and $g \in c$, where $\pi_S^H \colon V^H \to V^S$ is the projection map.
Observe that $\tau_H \colon V^H \to V^H$ is a linear cellular automaton with memory set $S$ over the group $H$.
We then have
\begin{equation}\label{e;tau-tau-c}
\tau = \prod_{c \in G/H}\tau_c
\end{equation}
in the sense that
$\pi_c \circ \tau = \tau_c \circ \pi_c$
for all $c \in G/H$.
 
From \eqref{e;tau-tau-c} we immediately deduce that $Y = \prod_{c \in G/H} Y_c$, where $Y_c = \pi_c(Y) = \tau_c(V^c)$,
and that $\tau_c$ is injective for all $c \in G/H$.

Consider the inverse map $(\tau_H)^{-1} \colon Y_H \to V^H$. By countability of $H$, we have by \cite[Theorem 3.1]{CC-inj}
that $(\tau_H)^{-1}$ is a linear cellular automaton. Let $T \subset H$ be a memory set for $(\tau_H)^{-1}$ and let $\nu \colon \pi_T^H(Y_H) = \pi_T(Y)  \to V$ be the corresponding local defining map, where $\pi_T^H \colon V^H \to V^T$ is the projection map. Consider the linear cellular automaton $\sigma \colon Y \to V^G$
defined by setting
$$
\sigma(y)(g) = \nu(\pi_T(g^{-1}y))
$$
for all $y \in Y$ and $g \in G$. Note that if $\sigma_c \colon Y_c \to V^c$ is the map defined by setting 
$\sigma_c(z)(g) = \nu(\pi_T^H(g^{-1}z))$ for all $z \in Y_c$ and $g \in c$ then $\sigma = \prod_{c \in G/H} \sigma_c$.
In particular, we have
\begin{equation}
\label{e;s-t}
\sigma_H = (\tau_H)^{-1}.
\end{equation}

For all $c \in G/H$ one has $\sigma_c = (\tau_c)^{-1}$.
Indeed, given $c \in G/H$, $z \in V^c$ and $g \in c$, we have
\[
\begin{split}
(\sigma_c \circ \tau_c)(z)(g) & = \sigma_c(\tau_c(z))(g)\\
& = g^{-1}\sigma_c(\tau_c(z))(1_G)\\
\mbox{(by \eqref{e;g-1-x-c})} \ & = \sigma_H(g^{-1}\tau_c(z))(1_G)\\
\mbox{(again by \eqref{e;g-1-x-c})} \ & = \sigma_H(\tau_H(g^{-1}z))(1_G)\\
\mbox{(by \eqref{e;s-t})} \ &  = (g^{-1}z)(1_G)\\
& = z(g).
\end{split}
\]
This shows that $\sigma_c \circ \tau_c$ equals the identity map ${\mbox {\rm Id}}_{V^c} \colon V^c \to V^c$.

We have proved that $\tau^{-1} = \sigma \colon Y \to V^G$ is a linear cellular automaton. In particular, $\tau^{-1}$ is
uniformly continuous by Proposition \ref{p;CA} and therefore $\tau$ is a uniform embedding.
\end{proof}
\par

\begin{proof}[Proof of Theorem \ref{t:L-surjunct-closed}]
Let $N \in \NN(\Gamma)$ and let $(N_i)_{i \in I}$ be a net in $\NN(\Gamma)$ converging to $N$.
Suppose that the groups $G_i = \Gamma/N_i$ are L$_\K$-surjunctive for all $i \in I$ and let us show that the group $G = \Gamma/N$ is also L$_\K$-surjunctive.
\par
Let $V$ be a finite-dimensional vector space over the field $\K$  and let $\tau \colon V^G \to V^ G$ be an injective linear cellular automaton over $V^G$.
Let $\rho^* \colon V^G \to \Fix(N)$ be the linear uniform isomorphism (cf. Proposition \ref{p:describes-fix-N-lin}) given by 
defined by $\rho^* = y \circ \rho$ for all $y \in V^G$, where $\rho \colon \Gamma \to G$ is the canonical epimorphism.
By Corollary \ref{c:rho-conj-lca}, the map $\tau^* \colon \Fix(N) \to \Fix(N)$ given by $\tau^* = \rho^* \circ \tau \circ (\rho^*)^{-1}$ is a linear cellular automaton over the linear subshift $\Fix(N) \subset V^\Gamma$. 
By applying Proposition \ref{p;LCA}, 
we deduce the existence of a linear cellular automaton $\widetilde{\tau} \colon V^\Gamma \to V^\Gamma$ such that $\widetilde{\tau}(\Fix(N)) \subset \Fix(N)$ and $\widetilde{\tau}\vert_{\Fix(N)} = \tau^*$.
\par
 We claim that the hypotheses of Corollary \ref{c;gil} are satisfied by taking $X = V^\Gamma$, $f = \widetilde{\tau}$, $Y = \Fix(N)$, and $Z_i = \Fix(N_i)$. 
The action of $\Gamma$ on $X$ is uniformly continuous and expansive by Proposition \ref{p;UCE}.
Also, the linear cellular automaton  $f \colon X \to X$ is uniformly continuous and $\Gamma$-equivariant by Proposition \ref{p;LCA}.
The space $Y = \Fix(N)$ is closed in $X$. On the other hand, $\tau(V^G)$ is closed in $V^G$ by virtue of Lemma \ref{l;closure}, so that
 $f(Y) = \widetilde{\tau}(\Fix(N)) = \tau^*(\Fix(N))$ is closed in $Y = \Fix(N)$ and therefore in $X$. As $\tau$ is a injective linear cellular automaton, it follows from Lemma \ref{l;inverse} that $\tau$ is a uniform embedding of $V^G$ into itself.
Since $\rho^* \colon V^G \to \Fix(N)$ is a linear uniform isomorphism by  Proposition \ref{p:describes-fix-N-lin},
we deduce that $f\vert_Y = \widetilde{\tau}\vert_{\Fix(N)}= \tau^* =  \rho^* \circ \tau \circ (\rho^*)^{-1}$ is a uniform embedding.
We have $f(Z_i) \subset Z_i$ since the linear cellular automaton $f$ is $\Gamma$-equivariant. Moreover, $Z_i$ is $\Gamma$-invariant by normality of the
subgroup $N_i \subset \Gamma$.  By Theorem \ref{t:psi-uc-n-p}, we have $Z_i \overset{H-B}{\to} Y$ 
since $N_i$ converges to $N$. For $i \in I$, let $\rho_i \colon \Gamma \to G_i = \Gamma/N_i$ denote the canonical epimorphism and let $\rho_i^* \colon V^{G_i} \to \Fix(N_i)$ be the uniform isomorphism given by Proposition \ref{p:describes-fix-N-lin}. It follows from Corollary \ref{c:rho-conj-lca} that the map $\tau_i = (\rho_i^*)^{-1} \circ \widetilde{\tau} \vert_{Z_i} \circ \rho_i^* \colon V^{G_i} \to V^{G_i}$ is
a linear cellular automaton and hence is surjunctive since $G_i$ is an L$_\K$-surjunctive group. This implies that
$f\vert_{Z_i} = \widetilde{\tau}_{\Fix(N)}$ is surjunctive as well. 
By applying Corollary \ref{c;gil}, we deduce that $f(Y) = Y$. Therefore $\tau$ is surjective.
This shows that the group $G$ is L$_\K$-surjunctive.
\end{proof}


\end{document}